\documentclass[10pt]{article}
\usepackage{amsmath}
\usepackage{latexsym}
\usepackage{amssymb}
\usepackage{amsfonts}
\usepackage{amsthm}
\usepackage{tikz}
\newcommand*\circled[1]{\tikz[baseline=(char.base)]{\node[shape=circle,draw,inner sep=2pt] (char){#1};}}
\title{A PEDESTRIAN APPROACH TO COSSERAT/MAXWELL/WEYL THEORY}
\author{J.-F. Pommaret \\ CERMICS, Ecole des Ponts ParisTech,\\ 6/8 Av. Blaise Pascal, 77455 Marne-la-Vall\'ee Cedex 02, France \\
E-mail: jean-francois.pommaret@wanadoo.fr, pommaret@cermics.enpc.fr \\
URL: http://cermics.enpc.fr/$\sim$pommaret/home.html }
\date{  }
\begin{document}
\maketitle

\noindent
{\bf INTRODUCTION} :   \\

\[   \begin{array}{rcccc}
   &   &  CARTAN   &   \longrightarrow   &  SPENCER  \\
      &   \nearrow  &   &   &   \\
 LIE   &   &   \updownarrow   &  ?  &   \updownarrow   \\
      &   \searrow  &  &  &   \\
   &   &   VESSIOT   &   \longrightarrow   &  JANET   
   \end{array}     \]
\vspace*{1mm}  \\

\noindent
{\bf 1) LIE GROUPS $\longrightarrow $ LIE PSEUDOGROUPS }:  \\

\noindent
$X$= manifold with local coordinates $(x^i),i=1,...,n=dim(X)$ \\
$G$= Lie group with local coordinates  $(a^{\tau}),\tau=1,...,p=dim(G)$\\
{\it Lie group action} : $ X\times G \longrightarrow X \times X : (x,a) \longrightarrow (x,y=ax=f(x,a))$  \\
{\it Affine transformations} : $(x,(a^1,a^2)) \longrightarrow (x, y=a^1x+a^2)$  \hspace{5mm} (1 dilatation + 1 translation)  \\
{\it Finite transformations} : \hspace{1cm}$y=f(x)$ \hspace{13mm} $\longrightarrow $ ${\partial}_{xx}f(x)=0$  \hspace{2cm} (order 2)  \\
{\it Infinitesimal transformations} : $y=x+t\xi(x)+... $ $\longrightarrow $ $ {\partial}_{xx}\xi (x)=0 $ \hspace{2cm} (order 2)  \\
\noindent
{\it Remark} : The case of projective transformations needs more work and leads to order $3$ (exercise).  \\

\noindent
{\bf 2) VESSIOT $\longrightarrow $ JANET }:  \\

\noindent
{\it Jet coordinates at order 2} : $y_2=(y^k,y^k_i,y^k_{ij}) $ and $y_q$ at order $q$. \\
{\it Sections at order 2} : The section $f_2(x)=(f^k(x), f^k_i(x),f^k_{ij}(x))$ may be different from the section $j_2(f)(x) = (f^k(x),{\partial}_if^k(x),{\partial}_{ij}f^k(x))$ and the Spencer operator will measure the difference by setting $Df_2(x)=({\partial}_if^k(x)-f^k_i(x),{\partial}_if^k_j(x)-f^k_{ij}(x))$.\\
More generally, $f_q$ may be different from $j_q(f)$ and we shall simply set $Df_{q+1}=j_1(f_q)-f_{q+1}$.\\
{\it Differential invariants} $ \longrightarrow$ {\it Lie form} : $\Phi(y_2)\equiv \frac{y_{xx}}{y_x}= \omega (x) $ \hspace{1cm}  \\
\[ \bar{y}=a^1y+a^2 \Longrightarrow {\bar{y}}_x=a^1y_x \Longrightarrow {\bar{y}}_{xx}=a^1y_{xx} \Longrightarrow \Phi({\bar{y}}_2)=\Phi(y_2) \]
{\it Lie operator} $\longrightarrow $ {\it Lie derivative} : ${\cal{D}}\xi\equiv {\cal{L}}(\xi)\omega={\partial}_{xx}\xi=\Omega $ \\
${\cal{D}}\xi=0, {\cal{D}}\eta=0 \Longrightarrow {\cal{D}}[\xi,\eta]=0$ with bracket $[\xi ,\eta]=\xi{\partial}_x\eta-\eta{\partial}_x\xi$  \\

\noindent
{\bf 3) CARTAN $\longrightarrow $ SPENCER}:  \\

\noindent
{\it Roughly} : Arbitrary jet coordinates for vectors are considered as new unknowns and lead to {\it new sections} ${\xi}_2(x)=(\xi(x), {\xi}_x(x),{\xi}_{xx}(x)=0)$ such that the number of linearly independent components is equal to the number of parameters. \\
{\it General Spencer operator} :  \hspace{2mm}  $D: {\xi}_2\longrightarrow ({\partial}_x\xi(x)-{\xi}_x(x), {\partial}_x{\xi}_x(x)-{\xi}_{xx}(x))$  \hspace{1cm} (order 1)  \\
{\it Restricted Spencer operator} : $D: {\xi}_2 \longrightarrow ({\partial}_x\xi(x)-{\xi}_x(x), {\partial}_x{\xi}_x(x))$ \hspace{24mm} (order 1)\\
One has the following commutative diagram where the numbers of arbitrary functions are circled:  \\

\[  \begin{array}{rccccccccr}
  & && & &0 & & 0 &  &  \\
   & & & & & \downarrow & & \downarrow & &  \\
   & 0 &\longrightarrow &\Theta &\stackrel{j_2}{\longrightarrow} & \circled{\:2\:} &\stackrel{D_1}{\longrightarrow} &\circled{\:2\:} &\longrightarrow 0 & \hspace{3mm} Spencer \hspace{2mm} sequence  \\
   & & & & & \downarrow & & \parallel  &  &  \\
    & 0 & \longrightarrow &\circled{\:1\:} &\stackrel{j_2}{\longrightarrow} & \circled{\:3\:} &\stackrel{D_1}{\longrightarrow} &\circled{\:2\:} & \longrightarrow 0&   \\
    & & & \parallel & & \hspace{2mm}\downarrow \Phi & & \downarrow & &  \\
    0\longrightarrow & \Theta & \longrightarrow & \circled{\:1\:} & \stackrel{{\cal{D}}}{\longrightarrow} & \circled{\:1\:} & \longrightarrow & 0 & & Janet \hspace{2mm} sequence \\
    & & & & & \downarrow & & & &  \\
     & & & & & 0 & & & & 
 \end{array}  \]

\noindent 
In this diagram, which only depends on the left commutative square, the operator $j_2:\xi(x)\rightarrow (\xi(x)=\xi(x),{\partial}_x\xi(x)={\xi}_x(x),{\partial}_{xx}\xi(x)={\xi}_{xx}(x))$ has compatibility conditions $D_1{\xi}_2=0$ induced by $D$ and the space of solutions $\Theta $ of ${\cal{D}}=\Phi\circ j_2$ is generated over the constants by the infinitesimal generators ${\theta}_1=x{\partial}_x$ (dilatation) and ${\theta}_2={\partial}_x$ (translation) of the action. \\

\noindent
{\bf 4) DUALITY $\longrightarrow $ FORMAL ADJOINT}: \\

\noindent
{\it Roughly} : Contrary to what happens in the Janet sequence, the formal adjoint of the Spencer operator brings as many dual equations as the number of parameters (1 translation + 1 dilatation). \\
\[   \sigma ({\partial}_x\xi-{\xi}_x)+\mu{\partial}_x{\xi}_x(x)=-[ ({\partial}_x\sigma) \xi + ({\partial}_x\mu +\sigma){\xi}_x] + {\partial}_x(\sigma \xi +\mu {\xi}_x)\]
\noindent
{\it Cosserat equations} :  \hspace{5mm}      $ {\partial}_x\sigma=f  \hspace{2mm}, \hspace{2mm} {\partial}_x\mu + \sigma = m  $      \hspace{2cm} (equivalent "momenta" ) \\
{\it Remark} : The case of projective transformations is similar but needs more work (exercise).  \\ 
 
\noindent
{\bf 5) CONFORMAL GROUP $\longrightarrow $ COSSERAT/MAXWELL/WEYL EQUATIONS}:  \\ 
 
 \noindent
 Applying the above techniques to the group of conformal transformations of space-time transforming the Minkowski metric up to a function factor (15 parameters = 4 translations + 3 space rotations + 3 Lorentz transformations + 1 dilatation + 4 elations) brings Cosserat equations exactly on equal footing with Weyl equations and thus with Maxwell equations. This result provides for the first time the group theoretical unification of finite elements in engineering sciences. However, the previous methods are still not known by the mathematical and mechanical communities for reasons that are largely not scientific at all. As a byproduct, we do not know other references on these topics.\\
 
 \noindent
 {\it Example} : If we restrict our study to the group of isometries of the euclidean metric $\omega$ in dimension $n\geq 2$, exhibiting the Janet and the Spencer sequences is not easy at all, even when $n=2$, because the corresponding Killing operator ${\cal{D}}\xi={\cal{L}}(\xi)\omega=\Omega$, involving the Lie derivative ${\cal{L}}$ and providing twice the so-called infinitesimal deformation tensor $\epsilon$ of continuum mechanics, is not involutive. In order to overcome this problem, one must differentiate once by considering also the Christoffel symbols $\gamma$ and add the operator ${\cal{L}}(\xi)\gamma=\Gamma $ with the well known Levi-Civita isomorphism $j_1(\omega)\simeq (\omega,\gamma)$. Introducing the bundle ${\wedge}^rT^*$ of completely skewsymmetric covariant tensors or $r$-forms and the exterior derivative $d$ with $d^2=d\circ d\equiv 0$, we have the Poincar\'{e} sequence:\\
 \[   {\wedge}^0T^* \stackrel{d}{\longrightarrow} {\wedge}^1T^* \stackrel{d}{\longrightarrow} {\wedge}^2T^* \stackrel{d}{\longrightarrow}  ...  \stackrel{d}{\longrightarrow} {\wedge}^nT^* \longrightarrow 0  \]
{\it For Lie groups of transformations, one can prove that the Spencer sequence is locally isomorphic to the tensor product of the Poincar\'{e} sequence by the Lie algebra of the underlying Lie group}. Hence, the bigger is the group involved, the bigger are the dimensions of the Spencer bundles, contrary to what happens in the Janet sequence where the first Janet bundle has only to do with differential invariants. \\

{\it Remark} : This rather philosophical comment, namely to replace the Janet sequence by the Spencer sequence, must be considered as the crucial key for understanding the work of the brothers E. and F. Cosserat in 1909, the best picture being that of two children playing at see-saw. \\

When $n=2$, one has 3 parameters (2 translations + 1 rotation) and the following commutative diagram which, as before, only depends on the left commutative square:  \\
  
  \[  \begin{array}{rccccccccccccr}
 &&&&& 0 &&0&&0&  & \\
 &&&&& \downarrow && \downarrow && \downarrow  &\\
  & 0& \longrightarrow& \Theta &\stackrel{j_2}{\longrightarrow}&\circled{\:3\:}&\stackrel{D_1}{\longrightarrow}&\circled{\:6\:} &\stackrel{D_2}{\longrightarrow} & \circled{\:3\:}&\longrightarrow  0 & \hspace{3mm}Spencer \hspace{2mm} sequence \\
  &&&&& \downarrow & & \downarrow & & \downarrow & &    & \\
   & 0 & \longrightarrow & \circled{\:2\:} & \stackrel{j_2}{\longrightarrow} & \circled{12}& \stackrel{D_1}{\longrightarrow} & \circled{16} &\stackrel{D_2}{\longrightarrow} & \circled{\:6\:} &   \longrightarrow 0 &\\
   & & & \parallel && \hspace{5mm}\downarrow {\Phi}_0 & &\hspace{5mm} \downarrow {\Phi}_1 & & \hspace{5mm}\downarrow {\Phi}_2 &  &\\
   0 \longrightarrow & \Theta &\longrightarrow & \circled{\:2\:} & \stackrel{\cal{D}}{\longrightarrow} & \circled{\:9\:} & \stackrel{{\cal{D}}_1}{\longrightarrow} & \circled{10} & \stackrel{{\cal{D}}_2}{\longrightarrow} & \circled{\:3\:}& \longrightarrow  0 & \hspace{7mm} Janet \hspace{2mm} sequence\\
   &&&&& \downarrow & & \downarrow & & \downarrow &      &\\
   &&&&& 0 && 0 && 0  &  &
   \end{array}     \]
   
More generally, for $n\geq 2$ arbitrary, the adjoint of the first Spencer operator $D_1$ provides the Cosserat equations which can be parametrized by the adjoint of the second Spencer operator $D_2$ because it is well known that the Poincar\'{e} sequence is self-adjoint up to sign. A delicate 
theorem of homological algebra on the vanishing of the so-called {\it extension modules} (See [2] for more details) finally proves that the adjoint of the Lie operator ${\cal{D}}$ (stress equations) can also be parametrized by the adjoint of its compatibility conditions ${\cal{D}}_1$ (Airy functions).  \\

\noindent
{\it Remark} :  It is important to notice that the parametrization of the Cosserat equatons is thus {\it first order} while the parametrization of the classical stress equations (Airy when $n=2$, Morera/Maxwell when $n=3$) is {\it second order}, a result not evident at all which does not seem to be known today by mechanicians. \\

When $n=2$, the adjoint of $D_1$ provides the Cosserat equations. Indeed, lowering the upper indices of ${\xi}_2$ by means of the constant euclidean metric, we just need to look for the factors of ${\xi}_1,{\xi}_2$ and ${\xi}_{1,2}$ in the integration by parts of the sum:\\
 \[ {\sigma}^{11}({\partial}_1{\xi}_1-{\xi}_{1,1})+{\sigma}^{12}({\partial}_2{\xi}_1-{\xi}_{1,2})+{\sigma}^{21}({\partial}_1{\xi}_2-{\xi}_{2,1})+{\sigma}^{22}({\partial}_2{\xi}_2-{\xi}_{2,2})+{\mu}^{r}({\partial}_r{\xi}_{1,2}-{\xi}_{1,2r}) \]
 in order to obtain: \\
 
 $ {\partial}_1{\sigma}^{11}+{\partial}_2{\sigma}^{12}=f^1, {\partial}_1{\sigma}^{21}+{\partial}_2{\sigma}^{22}=f^2, {\partial}_1{\mu}^1+{\partial}_2{\mu}^2+{\sigma}^{12}-{\sigma}^{21}=m$ \hspace{2mm} (equivalent "momenta") \\

 Finally, we obtain the nontrivial {\it first order} parametrization ${\sigma}^{11}={\partial}_2{\phi}^1, {\sigma}^{12}=-{\partial}_1{\phi}^1, {\sigma}^{21}=-{\partial}_2{\phi}^2, {\sigma}^{22}={\partial}_1{\phi}^2, {\mu}^{1}={\partial}_2{\phi}^3+{\phi}^1, {\mu}^{2}=-{\partial}_1{\phi}^3-{\phi}^2$ in a coherent way with the Airy {\it second order} parametrization obtained if we set  ${\phi}^1={\partial}_2{\phi}, {\phi}^2={\partial}_1{\phi}, {\phi}^3=-\phi$ when 
 ${\mu}^1=0,{\mu}^2=0$.\\

\noindent
Meanwhile, the adjoint of the {\it second order} operator ${\cal{D}}_1: \epsilon \longrightarrow {\partial}_{11}{\epsilon}_{22}+{\partial}_{22}{\epsilon}_{11}-2{\partial}_{12}{\epsilon}_{12}$ is nothing else than the {\it second order} parametrization ${\sigma}^{11}={\partial}_{22}\phi, {\sigma}^{22}={\partial}_{11}\phi, {\sigma}^{12}={\sigma}^{21}=-{\partial}_{12}\phi$ of the classical stress equations by means of the single Airy function $\phi$.\\

The situation is even more delicate for the conformal group because certain properties are only existing for $n\geq 4$ though we have already 
${\xi}^k_{rij}=0, \forall n\geq 3$. Accordingly, when $n=4$, among the components of the Spencer operator we have ${\partial}_i{\xi}^r_{rj}- {\xi}^r_{rij}={\partial}_i{\xi}^r_{rj}$ and thus ${\partial}_i{\xi}^r_{rj}-{\partial}_j{\xi}^r_{ri}=F_{ij}$. Such a result allows to recover the electromagnetic (EM) field and Maxwell equations by duality along the way proposed by Weyl but the use of the Spencer operator provides a possibility to exhibit a link with Cosserat equations. \\ 

\noindent
{\bf 6) GAUGE THEORY (GT)}: \\

\noindent
{\it Gauging procedure} : If $y=a(t)x+b(t)$ with $a(t)$ a time depending orthogonal matrix ({\it rotation}) and $b(t)$ a time depending vector ({\it translation}) describes the movement of a rigid body in ${\mathbb{R}}^3$, then the projection of the speed $v=\dot{a}(t)x+\dot{b}(t)$ in an orthogonal frame fixed in the body is $a^{-1}v=a^{-1}\dot{a}x+a^{-1}\dot{b}$ and the kinetic energy is a quadratic function of the $1$-forms $a^{-1}\dot{a}$ and $a^{-1}\dot{b}$.  \\

More generally, we may consider a map $a: X \longrightarrow G:x \longrightarrow a(x)$, introduce the tangent mapping $T(a):T=T(X) \longrightarrow T(G):dx \longrightarrow da=\frac{\partial a}{\partial x} dx$ and consider the family of left invariant $1$-forms $a^{-1}da=A=(A^{\tau}_i(x)dx^i)$ with value in the {\it Lie algebra} ${\cal{G}}=T_e(G)$, the tangent space of $G$ at the identity $e\in G$ with {\it structure constants} $c=(c^{\tau}_{\rho\sigma})$. Using local coordinates, we may introduce the $2$-forms ${\partial}_iA^{\tau}_j-{\partial}_jA^{\tau}_i-c^{\tau}_{\rho\sigma}A^{\rho}_iA^{\sigma}_j=F^{\tau}_{ij}$ with value in ${\cal{G}}$, simply denoted by $dA-[A,A]=F$, and we have $A=a^{-1}da \Leftrightarrow F=0$ by pulling back on $X$ the {\it Maurer-Cartan equations} on $G$. \\

In 1956, at the birth of GT, the above notations were coming from the EM potential $A$ and EM field $dA=F$ of relativistic Maxwell theory. 
Accordingly, $G=U(1)$ (unit circle in the complex plane)$\longrightarrow dim ({\cal{G}})=1$ was the {\it only possibility} to get {\it pure} $1$-form $A$ and $2$-form $F$ when $c=0$.  \\

On the contrary, in the conformal framework where $G$ {\it is acting on} $X$, the second order jets ({\it elations}) ${\xi}^k_{ij}={\delta}^k_ia_j+{\delta}^k_ja_i-{\omega}_{ij}{\omega}^{kr}a_r \Rightarrow {\xi}^r_{ri}=na_i$ behave like the $1$-form $a_i(x)dx^i$ and the corresponding part of the Spencer operator $D$ is a $1$-form with value in $1$-form, that is a $(1,1)$-covariant tensor providing the EM field as a $2$-form by skewsymmetrization. This result, namely to construct lagrangians on the image of the induced Spencer operator $D_1$, is thus {\it perfectly coherent} with rigid body dynamics, Cosserat elasticity and Maxwell theory but in {\it total contradiction} with GT because $U(1)$ is not acting on space-time and there is a {\it shift by one step} in the interpretation of the Poincar\'{e} sequence involved because the fields are now described by $1$-forms.\\

\noindent
{\it Gauging procedure revisited} : Finally, we may extend the action $y=f(x,a)$ to $y_q=j_q(f)(x,a)$ in order to eliminate the parameters when $q$ is large enough. In this case, we may set $f(x)=f(x,a(x))$ and $f_q(x)=j_q(f)(x,a(x))$ in order to obtain $a(x)=a=cst \Leftrightarrow f_q=j_q(f)$ because 
$Df_{q+1}=j_1(f_q)-f_{q+1}=\frac{\partial f_q(x,a(x))}{\partial a^{\tau}}{\partial}_ia^{\tau}(x)$ and the matrix involved has maximum rank $p$.  \\

\noindent
{\bf 7) GENERAL RELATIVITY (GR)}:  \\

The mathematical foundation of GR is {\it always} presented in textbooks without any reference at all to conformal geometry and we first prove that such an approach is not correct indeed. For this, we shall compare the {\it classical Killing system} ${\Omega}_{ij}\equiv ( {\cal{L}}(\xi)\omega)_{ij}\equiv {\omega}_{rj}{\partial}_i{\xi}^r+{\omega}_{ir}{\partial}_j{\xi}^r+{\xi}^r{\partial}_r{\omega}_{ij}=0$ to the {\it conformal Killing system} ${\hat{\Omega}}_{ij}\equiv {\hat{\omega}}_{rj}{\partial}_i{\xi}^r+{\hat{\omega}}_{ir}{\partial}_j{\xi}^r-\frac{2}{n}{\hat{\omega}}_{ij}{\partial}_r{\xi}^r+{\xi}^r{\partial}_r{\hat{\omega}}_{ij}=0$ obtained by introducing ${\hat{\omega}}_{ij}={\omega}_{ij}/\mid det(\omega)\mid ^{\frac{1}{n}}$ or , equivalently, by eliminating $A(x)$ in ${\cal{L}}(\xi)\omega=A(x)\omega$. \\

Counting the number of derivatives of the $\xi$ at order $1,2,3$ and the number of derivatives of the $\Omega$ at order $0,1,2$ respectively when $\omega$ is a constant metric with $det(\omega)\neq 0$ (for example the Minkowski metric when $n=4$) while looking at the ranks of the corresponding matrices, we obtain by difference the number of {\it compatibility conditions} (CC) for $\Omega$, namely none at order $0$, none at order $1$ and $n^2(n+1)^2/4-n^2(n+1)(n+2)/6=n^2(n^2-1)/12$ at order $2$ ({\it Riemann tensor}).\\

Proceeding in the same way for the conformal case, we get $n(n+1)/2(n(n+1)/2-1)-n^2(n+1)(n+2)/6=n(n+1)(n+2)(n-3)/12$ CC of order $2$ for the $\hat{\Omega}$ when $n\geq 3$ ({\it Weyl tensor}). Also, as the group of isometries is a subgroup of the group of conformal isometries, {\it the Riemann tensor projects onto the Weyl tensor and the kernel of this canonical projection is the Ricci tensor} with $n(n+1)/2$ components. As a byproduct, the Ricci tensor only depends on the "difference" existing between the clasical Killing system and the conformal Killing system, namely the $n$ second order jets ({\it elations} once more). However, {\it apart from a delicate diagram chasing}, there is no simple explanation of the fact that the Ricci tensor, thus obtained without contracting the indices as usual, may be embedded in the image of the Spencer operator made by $1$-forms with value in $1$-forms that we have already exhibited for describing EM. \\

\noindent
{\it Remark} :  It follows that the foundations of both GR and GT are not coherent with jet theory and must therefore be revisited within this new framework.  \\

\noindent
{\bf CONCLUSION}:   \\

\noindent
These new unavoidable methods based on the formal theory of systems of partial differential equations and Lie pseudogroups provide the common secret of the three following famous books [C], [M] and [W] published about at the same time at the beginning of the last century. Indeed, the Spencer operator can always be exhibited even if there is no group background and, when only constant sections are considered, one recovers exactly (up to sign) the operator introduced by Macaulay for studying {\it inverse systems}. This short notice can also be considered as an elementary summary of certain recent results presented in the references below. \\
 
 \noindent
[C]  E. and F. COSSERAT: Th\'{e}orie des Corps D\'{e}formables, Hermann, Paris, 1909.  \\
\noindent
[M]  F.S. MACAULAY: The Algebraic Theory of Modular Systems, Cambridge, 1916.  \\
\noindent
[W]  H. WEYL: Space, Time, Matter, Berlin, 1918 (1922, 1958; Dover, 1952). \\ 

\noindent
{\bf REFERENCES}:  \\

 \noindent
[1]  J.-F. POMMARET: Lie Pseudogroups and Mechanics, Gordon and Breach, New York, 1988.\\
 \noindent
[2]  J.-F. POMMARET: Partial Differential Control Theory, Kluwer, 2001.  \\
\noindent
[3]  J.-F. POMMARET: Parametrization of Cosserat Equations, Acta Mechanica, 215, 2010, 43-55.\\
 \noindent
[4]  J.-F. POMMARET: Macaulay Inverse Systems Revisited, Journal of Symbolic Computations, 46 (2011), 1049-1069.  \\
\noindent
[5]  J.-F. POMMARET: Spencer Operator and Applications: From Continuum Mechanics to Mathematical Physics, in "Continuum Mechanics-Progress in Fundamentals and Engineering Applications", Dr. Yong Gan (Ed.), ISBN: 978-953-51-0447--6, InTech, 2012, Available from: \\
http://www.intechopen.com/books/continuum-mechanics-progress-in-fundamentals-and-engineerin-applications/spencer-operator-and-applications-from-continuum-mechanics-to-mathematical-physics  \\

\end{document}